\documentclass{amsart}
\usepackage{amssymb}
\usepackage{amssymb}
\usepackage{epsfig,pb-diagram}
\usepackage[matrix, arrow, curve]{xy}

\title{On the Combinatorics of Crystal Graphs, II. The Crystal Commutor}

\author{Cristian Lenart}

\address{Department of Mathematics and Statistics, State University of New York at Albany, Albany, NY 12222}
\email{lenart@albany.edu}


\thanks{Cristian Lenart was supported by National Science Foundation 
grant DMS-0403029}


\numberwithin{equation}{section}

\theoremstyle{plain}
\newtheorem{theorem}{Theorem}[section]
\newtheorem{proposition}[theorem]{Proposition}
\newtheorem{lemma}[theorem]{Lemma}
\newtheorem{corollary}[theorem]{Corollary}
\newtheorem{conjecture}[theorem]{Conjecture}

\theoremstyle{definition}

\newtheorem{example}[theorem]{Example}

\theoremstyle{remark}
\newtheorem{remark}[theorem]{Remark}
\newtheorem{remarks}[theorem]{Remarks}

\def\w{{\rm wt}}
\def\e{\widetilde{e}}
\def\f{\widetilde{f}}

\def\Z{\mathbb{Z}}

\def\h{\mathfrak{h}}
\def\hR{\mathfrak{h}^*_\mathbb{R}}

\newcommand{\casetwo}[3]{\left\{ \begin{array}{ll} #1 &\mbox{if $#2$} \\ #3 &\mbox{otherwise}\,. \end{array} \right.}
\newcommand{\casetwoc}[3]{\left\{ \begin{array}{ll} #1 &\mbox{if $#2$} \\ #3 &\mbox{otherwise}\,, \end{array} \right.}
\newcommand{\casetwoex}[4]{\left\{ \begin{array}{ll} #1 &\mbox{if $#2$} \\ #3 &\mbox{if $#4$} \,. \end{array} \right.}
\newcommand{\casetwoexc}[4]{\left\{ \begin{array}{ll} #1 &\mbox{if $#2$} \\ #3 &\mbox{if $#4$} \,, \end{array} \right.}

\newcommand{\opi}{\overline{\pi}}
\newcommand{\opip}{\overline{\pi}'}
\newcommand{\ope}{\overline{p}}
\newcommand{\opep}{\overline{p}'}
\newcommand{\oqu}{\overline{q}}
\newcommand{\oqup}{\overline{q}'}

\begin{document}
\bibliographystyle{plain}

\begin{abstract}
We present an explicit combinatorial realization of the commutor in the category of crystals which was first studied by Henriques and Kamnitzer. Our realization is based on certain local moves defined by van Leeuwen. 
\end{abstract}

\maketitle



\section{Introduction}

We work in the category $\mathfrak g$-{\tt Crystals} of crystals corresponding to representations of complex semisimple Lie algebras $\mathfrak g$. It is well-known that this is a {\em monoidal category} with an associative tensor product (e.g., see \cite{hakccc}). The crystals $A\otimes B$ and $B\otimes A$ are isomorphic via maps called {\em commutors}. The map 
\[{\rm flip}\::\:A\otimes B\rightarrow  B\otimes A\,,\;\;\;\;\;(a,b)\mapsto (b,a)\] is not a commutor. Henriques and Kamnitzer \cite{hakccc}, based on an idea of Berenstein, defined a commutor $\sigma_{A,B}$ based on Lusztig's involution on a crystal. They also proved that the category  $\mathfrak g$-{\tt Crystals} with this commutor is a {\em  coboundary category} (cf. {\rm \cite{driqha}}). More recently, Kamnitzer and Tingley \cite{katdcc} proved that the action of $\sigma_{A,B}$ on the highest weight elements (which determines it) is given by {\em Kashiwara's involution} on the Verma crystal; as remarked in \cite{katdcc}, Kashiwara's involution can be realized in terms of {\em Mirkovi\'c-Vilonen polytopes} \cite{kamcss}. Both of the above constructions of $\sigma_{A,B}$ depend on some maps of crystals whose explicit construction is nontrivial. In this paper, we present an explicit realization of the commutor $\sigma_{A,B}$ based on some local moves defined by van Leeuwen \cite{leuajt}; the latter generalize the realization of Sch\"utzenberger's {\em jeu de taquin} for Young tableaux (e.g., see \cite{fulyt}) via Fomin's {\em growth diagrams}  \cite[Appendix 1]{staec2}. Our realization is proved to hold, in particular, for all simple Lie algebras $\mathfrak g$ with the exception of those of type $E_8$, $F_4$, and $G_2$; moreover, it is conjectured to hold for all semisimple Lie algebras.

\section{Background on Crystals}
\label{sec:notation}

In this section, we recall some background information on finite root systems and crystal graphs. 

\subsection{Root systems}\label{rootsyst}

Let $\mathfrak{g}$ be a complex semisimple Lie algebra with Cartan subalgebra $\h$. Let $\varPhi\subset \h^*$ be its finite {\it root system} with positive roots $\varPhi^+$, and let $\hR\subset \h^*$ be the real span of $\varPhi$. Let $r$ be the rank of $\h$, and let $\alpha_i\in\varPhi^+$ for $i\in I:=\{1,\ldots,r\}$ be the corresponding {\it simple roots}. The nondegenerate scalar product on $\hR$ induced by the Killing form is denoted by $\langle\,\cdot\,,\,\cdot\,\rangle$. Given a root $\alpha$, we have the corresponding {\it coroot\/} $\alpha^\vee := 2\alpha/\langle\alpha,\alpha\rangle$. The reflection corresponding to a root $\alpha$ is denoted, as usual, by $s_\alpha$, and we have the {\it simple reflections\/} $s_i := s_{\alpha_i}$ for $i\in I$.  Let   $W\subset\mathrm{Aut}(\hR)$ be the Weyl group of $\varPhi$. The {\em length} of an element $w\in W$ is denoted by $\ell(w)$. Let $w_\circ$ be the longest element of the Weyl group (i.e., $\ell(w_\circ)=|\varPhi^+|$). 

We denote the {\it weight lattice\/} by $\varLambda$, the subset of {\em dominant weights\/} by $\varLambda^+$, and the {\it fundamental weights\/} by $\varLambda_1,\dots,\varLambda_r$. For each dominant weight $\lambda$, there is a finite dimensional irreducible representation of $\mathfrak g$ with highest weight $\lambda$, which is denoted by $V_\lambda$. Let $\Z[\varLambda]$ be the group algebra of the weight lattice $\varLambda$, which has
a $\Z$-basis of formal exponents $\{e^\lambda \::\: \lambda\in\varLambda\}$, and 
multiplication $e^\lambda\cdot e^\mu := e^{\lambda+\mu}$; in other words,
$\Z[\varLambda]=\Z[e^{\pm\varLambda_1},\cdots,e^{\pm\varLambda_r}]$ is the algebra of
Laurent polynomials in  $r$ variables.  
The formal characters of the modules $V_{\lambda}$ 
are given by $ch(V_{\lambda})=\sum_{\mu\in\varLambda} 
m_{\lambda}(\mu)\,e^\mu\in\Z[\varLambda]$,
where $m_{\lambda}(\mu)$ is the multiplicity of the weight $\mu$ in 
$V_{\lambda}$.

A weight $\lambda$ is called {\em minuscule} if $\langle\lambda,\alpha^\vee\rangle\in\{0,\pm 1\}$ for all $\alpha\in\varPhi$, and {\em quasi-minuscule} if $\langle\lambda,\alpha^\vee\rangle\in\{0,\pm 1,\pm 2\}$ for all $\alpha\in\varPhi$. Since $W$ permutes $\varPhi$, every element in the $W$-orbit of a minuscule (or quasi-minuscule) weight is also minuscule (respectively quasi-minuscule). The irreducible representations corresponding to the dominant minuscule weights are called {\em minuscule representations}, and are characterized by the fact that their weights form a single $W$-orbit (i.e., the extremal weights, of multiplicity 1). The dominant minuscule weights have a well-known classification (e.g., see \cite[Exercise VI.4.15]{bougal}). 
There are minuscule weights for all Cartan-Killing types with the exception of types $E_8$, $F_4$, and $G_2$. 

We will need the following result in \cite{stecmw}.

\begin{theorem}\cite{stecmw}\label{decomp} If $\varPhi$ is irreducible, then for every dominant weight $\lambda$ there is a decomposition $\lambda=\lambda_1+\ldots+\lambda_n$ such that $\lambda_i$ is minuscule or quasi-minuscule, and $\lambda_1+\ldots+\lambda_i$ is dominant for all $i=1,\ldots,n$. In fact, all $\lambda_i$ can be chosen to be minuscule if $\varPhi$ is not of type $E_8$, $F_4$, and $G_2$. Furthermore, all $\lambda_i$ can be chosen in the orbit of a single minuscule or quasi-minuscule weight unless $\varPhi$ is of type $D_n$ for $n$ even.
\end{theorem}

\subsection{Crystals}
\label{sec:crystals}

This section follows \cite[Section 2]{stecmw} and \cite[Section 2]{hakccc}. We refer to these papers for more details.

Kashiwara \cite{kascqa,kascb} defined a directed colored graph, called {\em crystal graph}, on the {\em canonical basis} of a representation of the quantum group $U_q(\mathfrak g)$ \cite{kascqa,luscba}. This graph partially encodes the action of the Chevalley generators of $U_q(\mathfrak g)$ on the canonical basis.  We will now define axiomatically the category  $\mathfrak g$-{\tt Crystals} of crystals corresponding to complex semisimple Lie algebras $\mathfrak g$. 

The objects of the category $\mathfrak g$-{\tt Crystals}, which are called {\em crystals} (or ${\mathfrak g}$-{\em crystals}), are 4-tuples \linebreak  $(B,\w,\varepsilon,\{\f_1,\ldots,\f_r\})$, where
\begin{itemize}
\item $B$ is the underlying set of the corresponding crystal;
\item $\w$ and $\varepsilon$ are maps $B\rightarrow\varLambda$;
\item $\f_i$ are maps from $B$ to $B\sqcup \{0\}$.
\end{itemize}
For each $b\in B$, we call $\w(b)$, $\varepsilon(b)$, and $\phi(b):=\w(b)+\varepsilon(b)$ the {\em weight}, {\em depth}, and {\em rise} of $b$. A morphism of crystals is a map of the underlying sets that commutes with all the structure maps. 

By a result of Joseph \cite{josqgp}, the category $\mathfrak g$-{\tt Crystals} is uniquely defined by the axioms below. 

\vspace{2mm}
{\bf (A1)} $\varepsilon(b),\,\phi(b)\in\varLambda^+$.
\vspace{2mm}

We define the depth and rise in the direction of the simple root $\alpha_i$ by $\varepsilon_i(b):=\langle \varepsilon(b),\alpha_i^\vee\rangle$ and $\phi_i(b):=\langle \phi(b),\alpha_i^\vee\rangle$. We also let $\w_i(b):=\langle \w(b),\alpha_i^\vee\rangle$.

\vspace{2mm}
{\bf (A2)} $\f_i$ is a bijection from $\{b\in B\::\: \phi_i(b)>0\}$ to $\{b\in B\::\: \varepsilon_i(b)>0\}$.
\vspace{2mm}

We let $\e_i:=\f_i^{-1}$ denote the inverse map, and extend it to a map from $B$ to $B\sqcup\{0\}$ by defining it to be 0 on $\{b\in B\::\: \varepsilon_i(b)=0\}$.

\vspace{2mm}
{\bf (A3)} We have $\w(\f_i(b))=\w(b)-\alpha_i$, and $\varepsilon_i(\f_i(b))=\varepsilon_i(b)+1$, whenever $\f_i(b)\ne 0$.
\vspace{2mm}

Hence, we also have $\phi_i(\f_i(b))=\phi_i(b)-1$.  The maps $\e_i$ and $\f_i$, called {\em root operators}, act as raising and lowering operators which provide a partition of $B$ into $\alpha_i$-strings that are closed under the action of $\e_i$ and $\f_i$. For example, the $\alpha_i$-{\em string through} $b$ is (by definition)
\[\f_i^\phi(b),\,\ldots,\,\f_i(b),\,b,\,\e_i(b),\,\ldots,\,\e_i^{\varepsilon}(b)\,,\]
where $\varepsilon:=\varepsilon_i(b)$ and $\phi:=\phi_i(b)$; furthermore, we have $\f_i^{\phi+1}(b)=\e_i^{\varepsilon+1}(b)=0$. 

We define partial orders on $B$, one for each $i\in I$, by
\begin{equation}\label{partialord}c\preceq_ib\;\;\;\mathrm{if}\;\;\;c=\f_i^k(b)\;\:\textrm{for some }\,k\ge 0\,.\end{equation}
Let $\preceq$ denote the partial order on $B$ generated by all partial orders $\preceq_i$, for $i=1,\ldots,r$. The set of maximal elements of the poset $(B,\preceq)$ is denoted $\max\,B$. 

The {\em direct sum} of two crystals is defined by the disjoint union of their underlying sets. 

\vspace{2mm}
{\bf (A4)} For each dominant weight $\lambda$, the category contains an object $B_\lambda$ such that the poset $(B_\lambda,\preceq)$ has a maximum (i.e., a highest weight element), denoted $b_\lambda$, and $\w(b_\lambda)=\lambda$. Furthermore, the category consists of all crystals isomorphic to a direct sum of crystals $B_\lambda$. 
\vspace{2mm}

The crystal $B_\lambda$ is called a {\em highest weight crystal of highest weight} $\lambda$, and $b_\lambda$ is called its highest weight element. 

The {\em tensor product} of two crystals is the crystal corresponding to the tensor product of the corresponding representations; the explicit definition is given in \cite{kascb}. By iterating this construction, we can define the $n$-fold tensor product of crystals $A_1\otimes \ldots \otimes A_n$ as follows:
\begin{align}
&\w(a_1,\ldots,a_n):=\w(a_1)+\ldots+\w(a_n)\,,\nonumber\\
&\varepsilon_i(a_1,\ldots,a_n):=\max_{1\le j\le n}\varepsilon_i(a_j)-\w_i(a_1)-\ldots-\w_i(a_{j-1})\,,\label{defmd}\\
&\phi_i(a_1,\ldots,a_n):=\max_{1\le k\le n}\phi_i(a_k)+\w_i(a_{k+1})+\ldots+\w_i(a_{n})\,.\label{defme}
\end{align}
Furthermore, the root operators are defined by
\begin{align*}
&\e_i(a_1,\ldots,a_n):=(a_1,\ldots,\e_i(a_j),\ldots,a_n)\,,\\
&\f_i(a_1,\ldots,a_n):=(a_1,\ldots,\f_i(a_k),\ldots,a_n)\,,
\end{align*}
where $j$ and $k$ are the smallest and largest indices for which the maximum is achieved in (\ref{defmd}) (or, equivalently,  in (\ref{defme})). 

\vspace{2mm}
{\bf (A5)} For all dominant weights $\lambda,\mu$, there exists an inclusion of crystals $\iota_{\lambda,\mu}\::\:B_{\lambda+\mu}\hookrightarrow B_\lambda\otimes B_\mu$. 
\vspace{2mm}

The category $\mathfrak g$-{\tt Crystals} is closed under tensor products. Thus, we have
\begin{equation}\label{lrrule}
B_{\lambda_1}\otimes\ldots\otimes B_{\lambda_n}\simeq\bigoplus_{(b_1,\ldots,b_n))\in\max\, B_{\lambda_1}\otimes\ldots\otimes B_{\lambda_n}}B_{\w(b_1)+\ldots+\w(b_n)}\,,
\end{equation}
where $\simeq$ denotes isomorphism of crystals. By (\ref{defmd}), we have
\begin{equation}\label{maxchargen}
\max\,B_{\lambda_1}\otimes\ldots\otimes B_{\lambda_n}=\{(b_1,\ldots,b_n) \::\:\w(b_1)+\ldots+\w(b_{j-1})-\varepsilon(b_j)\in\varLambda^+\,,\;\:j=1,\ldots,n\}\,.
\end{equation}

\begin{remark}\label{lrmin} (1) If $\mu$ is a minuscule weight, then all the summands $B_\nu$ in $B_\lambda\otimes B_\mu$ correspond to $\nu=\lambda+\overline{\mu}$ for $\overline{\mu}\in W\mu$, and appear with multiplicity one.

(2) If $\mu$ is a quasi-minuscule weight, then all the summands $B_\nu$ in $B_\lambda\otimes B_\mu$ correspond to either $\nu=\lambda$, or $\nu=\lambda+\overline{\mu}$ for $\overline{\mu}\in W\mu$; the summands of the second type appear with multiplicity one.
\end{remark}

The above axioms imply that any highest weight crystal $B_\lambda$ has a minimum $c_\lambda$ with $\w(c_\lambda)=w_\circ(\lambda)$. Let $i\mapsto i^*$, for $i\in I$, be the Dynkin diagram automorphism specified by $\alpha_{i^*}=-w_\circ(\alpha_i)$. There is an involution $\eta_{B_\lambda}=\eta_\lambda$ on $B_\lambda$ specified by the following conditions:
\begin{equation}\label{invol}
\eta_\lambda(b_\lambda)=c_\lambda\,,\;\;\;\e_i\eta_\lambda=\eta_\lambda \f_i^*\,,
\end{equation}
where $\eta_\lambda(0)=0$. We can define an involution $\eta_B$ on any crystal $B$ by applying $\eta_\lambda$ to each component of $B$ isomorphic to $B_\lambda$. This involution is due to Lusztig, and corresponds to the action of $w_\circ$ on a representation of a quantum group (e.g., see \cite{bazcbq,lusiqg}). Thus, we have 
\begin{equation}\w(\eta_B(b))=w_\circ(\w(b))\,,\end{equation}
which also follows from the axioms. 

\subsection{Models for crystals}

There are several models for crystals corresponding to semisimple Lie algebras, such as: Kashiwara-Nakashima tableaux \cite{kancgr}, Littelmann paths \cite{litpro,litcrp}, the alcove path model \cite{lenccg,lapawg,lapcmc}, the model in \cite{baztpm} based on Lusztig's parametrization of canonical bases, as well as some models based on geometric constructions. Instead of using a particular model, we will use here the construction of crystals based on their embedding into tensor products of minuscule or quasi-minuscule crystals (i.e., crystals $B_\lambda$ with $\lambda$ minuscule or quasi-minuscule weights). We now describe this construction, following \cite{stecmw}. 

Let us first describe minuscule crystals. If $\lambda$ is minuscule, then the underlying set of $B_\lambda$ is the $W$-orbit $W\lambda$. Naturally, the weight of $\mu\in B_\lambda$ is $\mu$ itself, and we define
\begin{align*}
&\varepsilon_i(\mu):=\max(0,-\langle\mu,\alpha_i^\vee\rangle)=\casetwoc{1}{\langle\mu,\alpha_i^\vee\rangle=-1}{0}\\
&\phi_i(\mu):=\max(0,\langle\mu,\alpha_i^\vee\rangle)=\casetwo{1}{\langle\mu,\alpha_i^\vee\rangle=1}{0}
\end{align*}
The root operators are defined as follows, for any $\mu\in B_\lambda$:
\[\begin{array}{ll}
\f_i(\mu):=s_i(\mu)=\mu-\alpha_i\;\;\;\;&\mbox{if $\langle\mu,\alpha_i^\vee\rangle=1$}\,,\\[0.05in]
\e_i(\mu):=s_i(\mu)=\mu+\alpha_i\;\;\;\;&\mbox{if $\langle\mu,\alpha_i^\vee\rangle=-1$}\,.
\end{array}\]
In all other cases, the operators are defined to be $0$. 

For the description of a quasi-minuscule crystal, we refer to \cite{stecmw}. 

Now let $\lambda$ be any dominant weight. By Theorem \ref{decomp}, there is a decomposition $\lambda=\lambda_1+\ldots+\lambda_n$ such that $\lambda_i$ is minuscule or quasi-minuscule, and $\lambda_1+\ldots+\lambda_i$ is dominant for all $i=1,\ldots,n$. Let $\omega_i:={\rm dom}_W(\lambda_i)$ denote the dominant representative of the $W$-orbit $W\lambda_i$. Then, by Axiom (A5) and the construction of $n$-fold tensor products in Section \ref{sec:crystals}, we have an embedding
\begin{equation}\label{embed}B_\lambda\hookrightarrow B_{\omega_1}\otimes\ldots \otimes B_{\omega_n}\,.\end{equation}
Indeed, by (\ref{defmd}), the condition that all partial sums $\lambda_1+\ldots+\lambda_i$ are dominant weights guarantees that $(\lambda_1,\ldots,\lambda_n)$ is a maximal element in $B_{\omega_1}\otimes\ldots \otimes B_{\omega_n}$; thus, we can map $b_\lambda\mapsto (\lambda_1,\ldots,\lambda_n)$, and extend this to an embedding (\ref{embed}) by using lowering operators.

\subsection{A commutor for crystals}\label{cobcat}

This section follows \cite[Section 3]{hakccc}, and we refer to this paper for more details.

Henriques and Kamnitzer \cite{hakccc}, based on an idea of Berenstein, defined a commutor $\sigma_{A,B}$ (i.e., an isomorphism between the crystals $A\otimes B$ and $B\otimes A$) as follows:
\[\sigma_{A,B}(a,b):=\eta_{B\otimes A}(\eta_B(b),\,\eta_A(a))\,.\]
In other words, we have
\begin{equation}\label{defcomm}\sigma_{A,B}=\eta_{B\otimes A}\circ(\eta_B\otimes \eta_A)\circ\,{\rm flip}\,.\end{equation}
 It turns out that we also have 
\begin{equation}\label{newdefcomm}\sigma_{A,B}={\rm flip}\,\circ\eta_{A\otimes B}\circ(\eta_A\otimes\eta_B)\,.\end{equation}
More recently, Kamnitzer and Tingley \cite{katdcc} proved that the action of the above commutor on the highest weight elements (which determines it) is given by Kashiwara's involution on the Verma crystal \cite{kascb}. 

Henriques and Kamnitzer proved that the category  $\mathfrak g$-{\tt Crystals} with this commutor is a coboundary category (cf. {\rm \cite{driqha}}). This statement amounts to the following three properties:
\begin{enumerate}
\item[(C1)] $\sigma_{A,B}$ is an isomorphism of crystals, and is natural in  $A$ and $B$ (i.e., it commutes with the maps on $A\otimes B$ and $B\otimes A$ coming from maps of crystals $A\rightarrow C$ and $B\rightarrow D$);
\item[(C2)] $\sigma_{A,B}\circ\sigma_{B,A}=1$;
\item[(C3)] the following diagram commutes:
\begin{equation} \label{pentagon}
\begin{diagram}
\node[2]{A \otimes B \otimes C}\arrow[3]{e,t}{1 \otimes \sigma_{ B, C}}
\arrow{s,l}{\sigma_{A,B} \otimes 1}
\node[3]{A \otimes C \otimes B}\arrow{s,r}{\sigma_{A, C \otimes B}}\\
\node[2]{B \otimes A \otimes C}\arrow[3]{e,t}{\sigma_{B \otimes A, C}}
\node[3]{C \otimes B \otimes A}
\end{diagram}
\end{equation}
\end{enumerate}

Note that the third condition is  an analog of the hexagon axiom for braided categories. By analogy with braid groups acting on multiple tensor products in braided categories, there is a group, denoted by $J_n$ and called the $n$-{\em fruit cactus group}, which acts on $n$-fold tensor products in a coboundary category. This group is generated by $s_{p, q}$, for $1 \le p < q \le n$, subject to the following relations:
\begin{enumerate}
\item[(R1)] $s_{p,q}^2 = 1$;
\item[(R2)] $s_{p,q}  s_{k,l} = s_{k,l} s_{p,q}$ if $p < q$, $k < l$, and either $q<k$ or $l<p$;
\item[(R3)] $s_{p,q}  s_{k,l} = s_{i,j} s_{p,q}$ if $p\le k<l\le  q$, where $i = \widehat{s}_{p,q}(l)$, $j = \widehat{s}_{p,q}(k)$,
\end{enumerate}
and $\widehat{s}_{p,q}$ denotes the following involution in the symmetric group $S_n$: 
\begin{equation*}
\widehat{s}_{p,q} = \begin{pmatrix}
1 & \cdots & p -1 & p &\cdots& q & q+1 &\cdots& n \\ 
1 & \cdots& p-1 & q& \cdots &p & q+1 &\cdots &n   
\end{pmatrix}\,.
\end{equation*}

In the general context of coboundary categories, it is shown in \cite[Lemmas 3-4]{hakccc} that the group $J_n$ acts on $n$-fold tensor products $A_1\ldots A_n:=A_1\otimes\ldots A_n$ by letting its generators ``reverse intervals'':
\[s_{p,q}:A_1  \cdots A_{p-1}A_p A_{p+1} \cdots A_{q-1}A_q A_{q+1} \cdots  A_n \longrightarrow A_1  \cdots A_{p-1}A_q A_{q-1} \cdots A_{p+1}A_p A_{q+1} \cdots  A_n\,.\]
More precisely, based on the commutor, we first define natural isomorphisms denoted $\sigma_{p,r,q}$ for $1 \le p \le r < q \le n$ by:
\begin{gather*}
 (\sigma_{p,r,q})_{A_1, \dots, A_n} := 1 \otimes \sigma_{A_p \cdots A_r, A_{r+1} \cdots A_q} \otimes 1:\\ 
A_1 \cdots A_{p-1}A_p \cdots A_r A_{r+1} \cdots A_q A_{q+1}\cdots  A_n \longrightarrow
A_1 \cdots A_{p-1}A_{r+1} \cdots A_q A_p \cdots A_rA_{q+1} \cdots A_n.
\end{gather*}
We then define the action of $s_{p,q}$ recursively by $s_{p, p+1}:= \sigma_{p,p,p+1}$, and $s_{p,q}:= \sigma_{p,p,q} \circ s_{p+1, q}$ for $q-p > 1$. The following expression of $\sigma_{p,r,q}$ in terms of the generators of $J_n$ will be needed:
\begin{equation}\label{sigmasss}
\sigma_{p, r, q}=s_{p,q} \circ s_{r+1, q}  \circ s_{p,r}\,.
\end{equation}

\section{Van Leeuwen's Jeu de Taquin}\label{vljdt}

In this section, we translate the construction of van Leeuwen in \cite{leuajt} into the language of crystal graphs. We start with the following important remark related to commutors.

\begin{remark}\label{constrcomm}
In order to define a commutor, i.e., an isomorphism between the crystals $B_{\pi'}\otimes B_\pi$ and $B_\pi\otimes B_{\pi'}$ (for $\pi',\pi\in\varLambda^+$), it suffices to construct a weight-preserving bijection between $\max\,B_{\pi'}\otimes B_\pi$ and $\max\,B_\pi\otimes B_{\pi'}$. Clearly, such a bijection can be uniquely extended to a commutor via the action of the lowering operators. 
\end{remark}

A bijection of this type was constructed by van Leeuwen in \cite{leuajt}, and is described below using the above setup.

We start by embedding the crystals $B_{\pi'}$ and $B_{\pi}$ into tensor products of minuscule or quasi-minuscule crystals, as in (\ref{embed}):
\begin{equation}\label{embedpi}B_{\pi'}\hookrightarrow B_{\omega'_1}\otimes\ldots \otimes B_{\omega'_l}\,,\;\;\;\;\;\;B_\pi\hookrightarrow B_{\omega_1}\otimes\ldots \otimes B_{\omega_k}\,;\end{equation}
these embeddings are based on decompositions $\pi'=\pi'_1+\ldots+\pi'_l$ and $\pi=\pi_1+\ldots+\pi_k$, i.e., we have $\omega'_j:={\rm dom}_W(\pi'_j)$ and $\omega_i:={\rm dom}_W(\pi_i)$. 
From now on, we will identify the crystals $B_{\pi'}$ and $B_\pi$ with their embeddings in the corresponding tensor products. 

\begin{remark}\label{natembed} The above embedding is compatible with the commutor of Henriques and Kamnitzer due to the naturality of the latter, cf. its property (C1).
\end{remark}

Consider an element (cf. (\ref{maxchargen}))
\begin{equation}\label{lowerel}(b_{\pi'},p)=((\pi'_1,\ldots,\pi'_l),\,(p_1,\ldots,p_k))\in\max \,B_{\pi'}\otimes B_\pi\,,\end{equation}
which means that $p_i\in B_{\omega_i}$ for $i=1,\ldots,k$.
It will be mapped to an element 
\begin{equation}\label{upperel}(b_{\pi},p')=((\pi_1,\ldots,\pi_k),\,(p_1',\ldots,p_l'))\in\max \,B_\pi\otimes B_{\pi'}\,,\end{equation}
which means that $p_j'\in B_{\omega'_j}$ for $j=1,\ldots,l$. This map is weight-preserving, so we have 
\begin{equation}\label{weightpres}\pi'+\w(p_1)+\ldots+\w(p_k)=\pi+\w(p_1')+\ldots+\w(p_l')\,.
\end{equation}
 The idea is to use certain local moves which generalize the realization of Sch\"utzenberger's jeu de taquin for Young tableaux (e.g., see \cite{fulyt}) via Fomin's growth diagrams \cite[Appendix 1]{staec2}. 

We will define two matrices of elements in minuscule or quasi-minuscule crystals
\begin{align*}&h_{i,j}\in B_{\omega'_{j+1}}\,,\;\;\;\mbox{for $i=0,\ldots,k$ and $j=0,\ldots,l-1$}\,,\\
&v_{i,j}\in B_{\omega_{i+1}}\,,\;\;\;\:\mbox{for $i=0,\ldots,k-1$ and $j=0,\ldots,l$}\,.
\end{align*}
They are related by 
\begin{equation}\label{weightpresloc}
\w(h_{i,j})+\w(v_{i,j+1})=\w(v_{i,j})+\w(h_{i+1,j})\,,\;\;\;\mbox{for $i=0,\ldots,k-1$ and $j=0,\ldots,l-1$}\,.
\end{equation}
These elements should be thought of as labels for the horizontal and vertical segments joining successive points in a square lattice which lie in the interior or on the rectangle with vertices $(0,0)$ and $(k,l)$. Here we use the matrix notation, with $k$ indicating the row, and $l$ the column. More precisely, $h_{i,j}$ and $v_{i,j}$ label the horizontal and vertical segments for which $(i,j)$ is the left endpoint, respectively the bottom endpoint. Condition (\ref{weightpresloc}) guarantees that the sum of weights of elements on each minimum length path between $(0,0)$ and $(i,j)$ is the same, where $0\le i\le k$ and $0\le j\le l$. Let $\lambda^{[i,j]}$ be this sum (by default, we set $\lambda^{[0,0]}:=0$). We also require the following conditions:
\begin{align}\label{domcond1}
&(b_{\lambda^{[i,j]}},h_{i,j})\in\max\,\,B_{\lambda^{[i,j]}}\otimes B_{\omega'_{j+1}}\,,\;\;\;\mbox{for $i=0,\ldots,k$ and $j=0,\ldots,l-1$}\,,\\
&(b_{\lambda^{[i,j]}},v_{i,j})\in\max\,\,B_{\lambda^{[i,j]}}\otimes B_{\omega_{i+1}}\,,\;\;\;\:\mbox{for $i=0,\ldots,k-1$ and $j=0,\ldots,l$}\,.\label{domcond2}
\end{align}

We start by setting
\begin{equation}\label{incond}h_{0,j}:=\pi'_{j+1}\,,\;\;\;\;\mbox{for $j=0,\ldots,l-1$}\,,\;\;\;\;\mbox{and}\;\;\;\;v_{i,l}:=p_{i+1}\,,\;\;\;\;\mbox{for $i=0,\ldots,k-1$}\,.\end{equation}
The other elements $h_{i,j}$ and $v_{i,j}$ are defined by local moves, which are bijections
\[\max\,B_{\lambda^{[i,j]}}\otimes B_{\omega'_{j+1}}\otimes B_{\omega_{i+1}}\longleftrightarrow \max\,B_{\lambda^{[i,j]}} \otimes B_{\omega_{i+1}} \otimes B_{\omega'_{j+1}}\,,\]
for $i=0,\ldots,k-1$ and $j=0,\ldots,l-1$. More precisely, we have
\begin{equation}\label{locmove}
(b_{\lambda^{[i,j]}},h_{i,j},v_{i,j+1})\mapsto (b_{\lambda^{[i,j]}},v_{i,j},h_{i+1,j})\,.
\end{equation}
The procedure ends by reading off the elements in (\ref{upperel}) as follows:
\begin{equation}\label{final}p_{j+1}':=h_{k,j}\,,\;\;\;\;\mbox{for $j=0,\ldots,l-1$}\,.\end{equation}

We now describe the local moves (\ref{locmove}) as weight-preserving bijections
\begin{equation}\label{domlocal}\max\,B_{\kappa}\otimes B_{\omega'}\otimes B_{\omega}\longleftrightarrow \max\,B_{\kappa} \otimes B_{\omega} \otimes B_{\omega'}\,,\end{equation}
for $\kappa,\omega',\omega\in\varLambda^+$ and $\omega',\omega$ minuscule or quasi-minuscule. Such a bijection is given by 
\[(b_\kappa,h,v)\mapsto (b_\kappa,v',h')\,.\]
 Let 
\[\lambda:=\kappa+\w(h)\,,\;\;\;\;\nu:=\lambda+\w(v)\,,\;\;\;\;\mu:=\kappa+\w(v')\,.\]

When $\omega'$ and $\omega$ are both dominant minuscule weights, the local moves are easily described by the condition
\begin{equation}\label{defmoves}
\mu={\rm dom}_W(\kappa+\nu-\lambda)\,,
\end{equation}
where ${\rm dom}_W(\cdot)$ denotes, as above, the dominant representative in the corresponding $W$-orbit. In other words, after computing $\mu$ based on $h$ and $v$, we set
\[v':=\mu-\kappa\;\;\;\;\mbox{and}\;\;\;\;h':=\nu-\mu\,.\]
It is not hard to check (see \cite{leuajt}) that these moves reduce to Fomin's local moves for his growth diagrams \cite[Appendix 1]{staec2} in the case of the root system $A_{n-1}$, as long as one uses embeddings of crystals into tensor powers of $B_{\varepsilon_1}$ -- the crystal of the vector representation of $\mathfrak{sl}_n$. For the definition of the local moves in the case when at least one of the weights $\omega$ and $\omega'$ is quasi-minuscule, we refer to \cite{leuajt}. 

In \cite{leuajt}, it is proved that these moves are weight-preserving bijections between the sets in (\ref{domlocal}). Hence, conditions (\ref{weightpresloc}), (\ref{domcond1}), and (\ref{domcond2}) are satisfied. It is also proved that the above moves are reversible, that is, the inverse bijection in (\ref{domlocal}) is given by the same moves. Finally, the following result is proved, based on the initial conditions (\ref{incond}), cf. also (\ref{lowerel}):
\begin{align}\label{rootops1}
&(h_{i,0},\ldots,h_{i,l-1})\succeq (h_{i+1,0},\ldots,h_{i+1,l-1}) \;\;\;\;\;\mbox{in $B_{\pi'}$, for $i=0,\ldots,k-1$}\,,\\
&(v_{0,j+1},\ldots,v_{k-1,j+1})\preceq (v_{0,j},\ldots,v_{k-1,j}) \;\;\;\;\mbox{in $B_\pi$, for $j=0,\ldots,l-1$}\,;\label{rootops2}
\end{align}
here we use the partial order on crystals defined in Section \ref{sec:crystals}. The above results imply (cf. also (\ref{maxchargen})):
\[((v_{0,0},\ldots,v_{k-1,0}),(h_{k,0},\ldots,h_{k,l-1}))=((\pi_1,\ldots,\pi_k),(h_{k,0},\ldots,h_{k,l-1}))\in\max\, B_{\pi}\otimes B_{\pi'}\,.\]
So (\ref{upperel}) is true based on (\ref{final}), while (\ref{weightpres}) follows from (\ref{weightpresloc}). The reversibility of the local moves implies that the inverse map, from $\max\,B_\pi\otimes B_{\pi'}$ to $\max\,B_{\pi'}\otimes B_\pi$, is given by the same local moves.

\section{The Explicit Realization of the Commutor}

In this section, we use the setup in Section \ref{vljdt} without further comment. Assume that we have $(b_{\pi'},p)\mapsto (b_\pi,p')$ by van Leeuwen's jeu de taquin (where $p\in B_\pi$ and $p'\in B_{\pi'}$). As discussed above, we also have $(b_{\pi},p')\mapsto (b_{\pi'},p)$.

\begin{proposition}\label{casemin}
We have $\sigma_{B_{\pi'},B_{\pi}}(b_{\pi'},p)=(b_{\pi},p')$ and $\sigma_{B_{\pi},B_{\pi'}}(b_{\pi},p')=(b_{\pi'},p)$ in the following two cases, which correspond to the special case $l=1$: {\rm (1)} $\pi'$ is a minuscule weight; {\rm (2)} $\pi'$ is a quasi-minuscule weight and $p'\in W\pi'$.
\end{proposition}

\begin{proof}
This is immediate by Remark \ref{lrmin}. Indeed, in both cases there is a unique choice for the commutor.
\end{proof}

Propositions \ref{casemin} allows us to express the local moves of van Leeuwen (\ref{domlocal}) in terms of the commutor of Henriques and Kamnitzer. We also need to recall the naturality of the latter, namely its property (C1).

\begin{corollary}\label{vlsigma} If $\omega'$ and $\omega$ are minuscule weights, the local move of van Leeuwen from $\max\,B_{\kappa}\otimes B_{\omega'}\otimes B_{\omega}$ to $\max\,B_{\kappa} \otimes B_{\omega} \otimes B_{\omega'}$ is given by the following composite:
\[ \sigma_{B_{\omega'},B_{\kappa}\otimes B_{\omega}}\circ(\sigma_{B_{\kappa},B_{\omega'}}\otimes 1)\,.\]
\end{corollary}

The following lemma is based on the notation related to coboundary categories that was introduced in Section \ref{cobcat}.

\begin{lemma}\label{cobcomm} The following two composites in a coboundary category coincide:
\begin{equation}
\begin{diagram}
\node{A_1A_2A_3A_4}\arrow[2]{e,t}{\sigma_{1,1,2}}\node[2]{A_2A_1A_3A_4}\arrow[2]{e,t}{\sigma_{1,1,3}}\node[2]{A_1A_3A_2A_4}\arrow[2]{e,t}{\sigma_{1,2,4}}\node[2]{A_2A_4A_1A_3}\arrow[2]{e,t}{\sigma_{1,2,3}}\node[2]{A_1A_2A_4A_3\,,}
\end{diagram}
\end{equation}
\begin{equation}
\begin{diagram}
\node{A_1A_2A_3A_4}\arrow[2]{e,t}{\sigma_{1,2,3}}\node[2]{A_3A_1A_2A_4}\arrow[2]{e,t}{\sigma_{1,1,4}}\node[2]{A_1A_2A_4A_3\,.}
\end{diagram}
\end{equation}
\end{lemma}

\begin{proof} By (\ref{sigmasss}), we have
\[\sigma_{1,1,3}\circ\sigma_{1,1,2}=s_{1,3}\circ s_{2,3}\circ s_{1,2}\,,\;\;\;\;\;\;\;\;\sigma_{1,1,2}\circ\sigma_{1,2,3}=s_{1,2}\circ s_{1,3}\circ s_{1,2}\,.\]
But, by relation (R3) in the cactus group (see Section \ref{cobcat}), we have $s_{1,3}\circ s_{2,3}=s_{1,2}\circ s_{1,3}$, which implies 
\begin{equation}\label{cob1}\sigma_{1,1,3}\circ\sigma_{1,1,2}=\sigma_{1,1,2}\circ\sigma_{1,2,3}\,.\end{equation}
Upon viewing the tensor product $A_3A_4$ as a single object of the category, the above relation also implies the following one (both sides being defined on $A_1A_2A_3A_4$):
\begin{equation}\label{cob2}\sigma_{1,1,4}\circ\sigma_{1,1,2}=\sigma_{1,2,3}\circ\sigma_{1,2,4}\,.\end{equation}
Finally, by (\ref{cob1}) and (\ref{cob2}), we have
\[\sigma_{1,2,3}\circ\sigma_{1,2,4}\circ\sigma_{1,1,3}\circ\sigma_{1,1,2}=(\sigma_{1,1,4}\circ\sigma_{1,1,2})\circ(\sigma_{1,1,2}\circ\sigma_{1,2,3})=\sigma_{1,1,4}\circ\sigma_{1,2,3}\,.\]
\end{proof}

We now arrive at our main result, which we state and prove in the setup of Section \ref{vljdt}. At this point, let us recall Remark \ref{constrcomm} related to the fact that van Leeuwen's jeu de taquin can be viewed as a commutor. Let us also recall Remark \ref{natembed} related to the embedding of a crystal into a tensor product of minuscule or quasi-minuscule crystals, which will be used without further comment in the proof below. 

\begin{theorem}\label{comagree}
Assume that all weights $\omega_j'$, for $j=1,\ldots,l$, and $\omega_i$, for $i=1,\ldots,k$, are minuscule. We have
\[\sigma_{B_{\pi'},B_{\pi}}(b_{\pi'},p)=(b_{\pi},p')\,,\;\;\;\;\;\mbox{and}\;\;\;\;\;\sigma_{B_{\pi},B_{\pi'}}(b_{\pi},p')=(b_{\pi'},p)\,.\]
\end{theorem}

\begin{proof}
We proceed by induction on $k+l$, also assuming that the result holds for $l=1$ or $k=1$, by Propositions \ref{casemin}.

Let 
\[\opip:=\pi_1'+\ldots+\pi_{l-1}'\,,\;\;\;\;\;\;\;\;\opi:=\pi_1+\ldots+\pi_{k-1}\,.\]
Let us realize the maxima in the crystals $B_{\opip}$ and $B_{\opi}$ as
\[b_{\opip}=(\pi_1',\ldots,\pi_{l-1}')\,,\;\;\;\;\;\;\;\;b_{\opi}=(\pi_1,\ldots,\pi_{k-1})\,.\]
Furthermore, let us define
\[\ope:=(p_1,\ldots,p_{k-1})\,,\;\;\;\;\;\;\;\;\opep:=(p_1',\ldots,p_{l-1}')\,.\]
Let us also consider the elements $q\in B_{\pi}$ and $q'\in B_{\pi'}$ (cf. (\ref{rootops1}) and (\ref{rootops2})) defined as follows:
\begin{align*}&q=(q_1,\ldots,q_k)\;\;\;\;\;\mbox{for $q_{i+1}:=v_{i,l-1}$}\,,\;\;\;\;\;\mbox{where $i=0,\ldots,k-1$, and}\\&q':=(q_1',\ldots q_{l}')\;\;\;\;\;\,\mbox{for $q_{j+1}':=h_{k-1,j}$}\,,\;\;\;\mbox{where $j=0,\ldots,l-1$.}\end{align*}
In a similar way to the definition of $\ope$ and $\opep$, we define $\oqu$ and $\oqup$. Let us now express the van Leeuwen map $(b_{\pi'},p)\mapsto (b_{\pi},p')$ as the following composite
\begin{align}\label{vlcomp}
&(b_{\pi'},p)=(b_{\opip},\pi_l',\ope,p_k)\mapsto (b_{\opi},\oqup,q_l',p_k)\mapsto (b_{\opip},\oqu,q_l',p_k)\mapsto \\
\mapsto &(b_{\opip},\oqu,q_k,p_l')\mapsto (b_{\opi},\pi_k,\opep,p_l')=(b_{\pi},p')\,.\nonumber
\end{align}

Now let us denote
\[B_1:=B_{\omega_1'}\otimes\ldots\otimes B_{\omega_{l-1}'}\,,\;\;\;\;\;B_2:=B_{\omega_l'}\,,\;\;\;\;\;B_3:=B_{\omega_1}\otimes\ldots\otimes B_{\omega_{k-1}}\,,\;\;\;\;\;B_4:=B_{\omega_k}\,.\]
By induction, we know that the composition of the first two maps in (\ref{vlcomp}) is given by the commutor of Henriques and Kamnitzer, as follows (see the notation in Section \ref{cobcat}):
\begin{equation}\label{com1}
\begin{diagram}
\node{B_1B_2B_3B_4}\arrow[2]{e,t}{\sigma_{1,2,3}}\node[2]{B_3B_1B_2B_4}\arrow[2]{e,t}{\sigma_{1,1,2}}\node[2]{B_1B_3B_2B_4\,.}
\end{diagram}
\end{equation}
Similarly, the last map in (\ref{vlcomp}) is given by 
\begin{equation}\label{com2}
\begin{diagram}
\node{B_1B_3B_4B_2}\arrow[2]{e,t}{\sigma_{1,1,3}}\node[2]{B_3B_4B_1B_2\,.}\end{diagram}
\end{equation}
Finally, the third map in (\ref{vlcomp}) is just van Leeuwen's local move (\ref{domlocal}); by Corollary \ref{vlsigma}, it can be expressed as
\begin{equation}\label{com3}
\begin{diagram}
\node{B_1B_3B_2B_4}\arrow[2]{e,t}{\sigma_{1,2,3}}\node[2]{B_2B_1B_3B_4}\arrow[2]{e,t}{\sigma_{1,1,4}}\node[2]{B_1B_3B_4B_2\,.}
\end{diagram}
\end{equation}
  
Combining (\ref{com1}), (\ref{com2}), and (\ref{com3}), we can express the composite (\ref{vlcomp}) as the following composition in terms of the commutor of Henriques and Kamnitzer:
\[\sigma_{1,1,3}\circ(\sigma_{1,1,4}\circ\sigma_{1,2,3})\circ\sigma_{1,1,2}\circ\sigma_{1,2,3}= 
\sigma_{1,1,3}\circ(\sigma_{1,2,3}\circ\sigma_{1,2,4}\circ\sigma_{1,1,3}\circ\sigma_{1,1,2})\circ\sigma_{1,1,2}\circ\sigma_{1,2,3}=\sigma_{1,2,4}\,.\]
The first equality follows by rewriting the bracket based on Lemma \ref{cobcomm}, while the second one uses the fact that $\sigma_{1,1,3}=\sigma_{1,2,3}^{-1}$ (cf. property (C2) of the commutor, stated in Section \ref{cobcat}).   

\end{proof}

\begin{remarks}\label{genstat} (1) Any crystal can be embedded into a tensor product of minuscule crystals in the case of an irreducible root system different from the ones of type $E_8$, $F_4$, and $G_2$ (cf. Theorem \ref{decomp}). Hence, the combinatorial realization (in Theorem \ref{comagree}) of the commutor due to Henriques and Kamnitzer based on the local moves of van Leeuwen holds in all the Lie types previously mentioned. 

(2) Our construction is more efficient than the original definition of the commutor (\ref{defcomm}). Indeed, the latter uses Lusztig's involution three times. This involution on irreducible crystals $B_\lambda$ can be realized combinatorially (via {\em Sch\"utzenberger's evacuation} on Young tableaux in type $A$ \cite{bazcbq,fulyt}, and via a similar procedure based on the alcove path model in arbitrary type \cite{lenccg}). However, for a tensor product $B_{\pi'}\otimes B_\pi$ of two irreducible crystals, there is no nice realization of this involution other than the definition (\ref{invol}); see Example \ref{excomm} for the concrete way to realize it on a minimal element of $B_{\pi}\otimes B_{\pi'}$ (this is always the case if we compute the action of $\sigma_{B_{\pi'}\otimes B_\pi}$ on $\max\,B_{\pi'}\otimes B_\pi$). Also note that the realization of the commutor mentioned in \cite{katdcc} involves the Mirkovi\'c-Vilonen polytopes, which are more complex objects. 

(3) We can combine our main result with both realizations of the commutor of Henriques and Kamnitzer (in \cite{hakccc} and \cite{katdcc}) in order to obtain explicit constructions related to: (i) Lusztig's involution on highest weight elements in a tensor product of two irreducible crystals; (ii) Kashiwara's involution.
\end{remarks}

Remark \ref{genstat} (1) leads us to the following conjecture, which was checked experimentally.

\begin{conjecture} Theorem {\rm \ref{comagree}} holds without the restriction on the embeddings {\rm (\ref{embedpi})} of the crystals $B_{\pi'}$ and $B_{\pi}$, i.e., it holds for embeddings into tensor products of minuscule and quasi-minuscule crystals. 
\end{conjecture}

\begin{example}\label{excomm} Consider the root system $A_2$, $\pi'=2\varepsilon_1+\varepsilon_2$, and $\pi=2\varepsilon_1$. We use the standard embeddings $B_{\pi'}\hookrightarrow B_{\varepsilon_1}^{\otimes 3}$ and $B_{\pi}\hookrightarrow B_{\varepsilon_1}^{\otimes 2}$; to be more precise,  $b_{\pi'}\mapsto(\varepsilon_1,\varepsilon_1,\varepsilon_2)$. Let
\[(b_{\pi'},p)=((\varepsilon_1,\varepsilon_1,\varepsilon_2),(\varepsilon_3,\varepsilon_2))\in\max\,B_{\pi'}\otimes B_\pi\,.\]
We claim that
\[\sigma_{B_{\pi'},B_{\pi}}(b_{\pi'},p)=((\varepsilon_1,\varepsilon_1),(\varepsilon_2,\varepsilon_2,\varepsilon_3))\in\max\,B_{\pi}\otimes B_{\pi'}\,.\]
Indeed, by Sch\"utzenberger's evacuation, we have $\eta_{B_{\pi'}}(b_{\pi'})=(\varepsilon_3,\varepsilon_2,\varepsilon_3)$, and $\eta_{B_{\pi}}(p)=(\varepsilon_2,\varepsilon_1)$. Since the element $((\varepsilon_2,\varepsilon_1), (\varepsilon_3,\varepsilon_2,\varepsilon_3))$ is a minimal element of $B_{\pi}\otimes B_{\pi'}$, by applying raising operators to it as long as possible, we obtain $((\varepsilon_1,\varepsilon_1),(\varepsilon_2,\varepsilon_2,\varepsilon_3))$. Alternatively, we can use the algorithm based on van Leeuwen's moves, which is shown in the diagram below. For simplicity, $\varepsilon_k$ is represented here by $k$; the labels on the horizontal segments are $h_{i\cdot}$, while those on the vertical ones are $v_{\cdot j}$, as indicated. Note that the algorithm inputs the labels on the bottom and right edges of the rectangle; it outputs those on the left and top edges, after labeling the segments in the interior. 
\[
    \epsfxsize=2.4in \epsffile{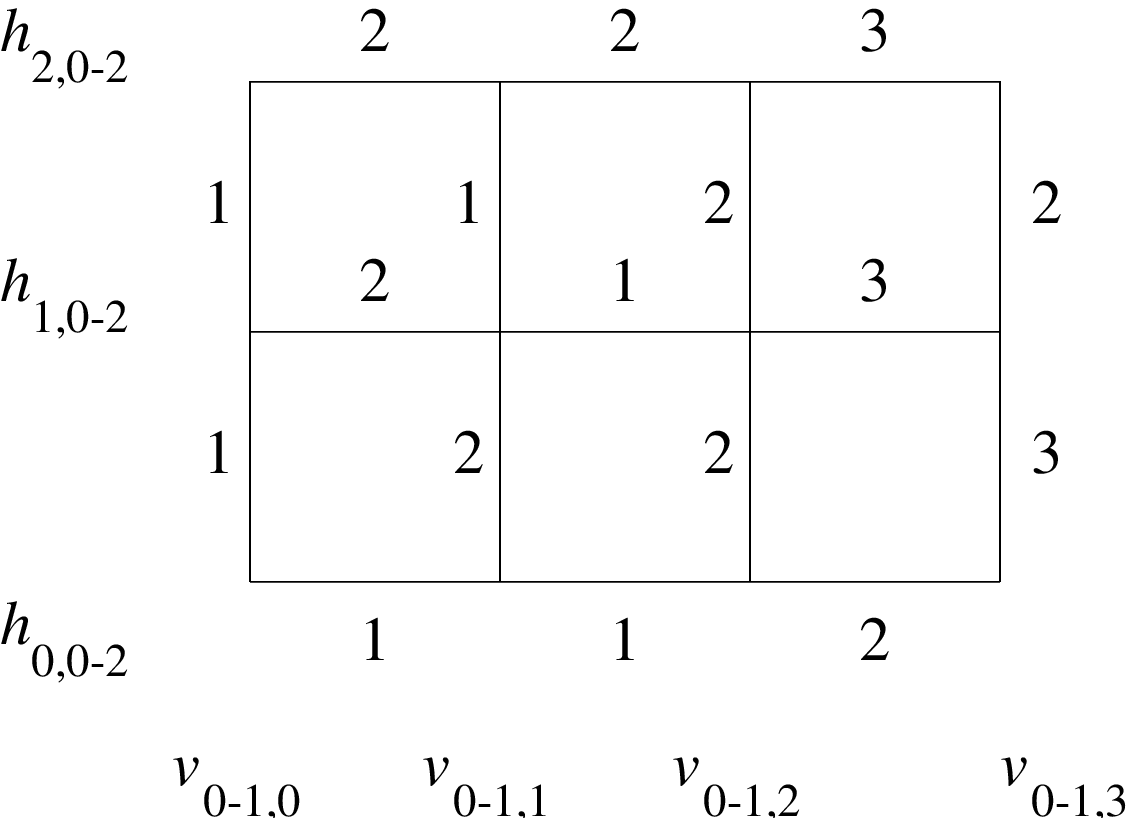}
\]
\end{example}


\begin{thebibliography}{10}

\bibitem{bazcbq}
A.~Berenstein and A.~Zelevinsky.
\newblock Canonical bases for the quantum group of type {$A\sb r$} and
  piecewise-linear combinatorics.
\newblock {\em Duke Math. J.}, 82:473--502, 1996.

\bibitem{baztpm}
A.~Berenstein and A.~Zelevinsky.
\newblock Tensor product multiplicities, canonical bases and totally positive
  varieties.
\newblock {\em Invent. Math.}, 143:77--128, 2001.

\bibitem{bougal}
N.~Bourbaki.
\newblock {\em Groupes et {A}lg{\`e}bres de {L}ie. \mbox{\rm {C}hp.
  {IV}--{VI}}}.
\newblock Masson, Paris, 1981.

\bibitem{driqha}
V.~G. Drinfel'd.
\newblock Quasi-{H}opf algebras.
\newblock {\em Leningrad Math. J.}, 1:1419--1457, 1990.

\bibitem{fulyt}
W.~Fulton.
\newblock {\em Young Tableaux}, volume~35 of {\em London Math. Soc. Student
  Texts}.
\newblock Cambridge Univ. Press, Cambridge and New York, 1997.

\bibitem{hakccc}
A.~Henriques and J.~Kamnitzer.
\newblock Crystals and coboundary categories.
\newblock {\em Duke Math. J.}, 132:191--216, 2006.

\bibitem{josqgp}
A.~Joseph.
\newblock {\em Quantum {G}roups and {T}heir {P}rimitive {I}deals}, volume~29 of
  {\em Ergebnisse der Mathematik und ihrer Grenzgebiete (3) [Results in
  Mathematics and Related Areas (3)]}.
\newblock Springer-Verlag, Berlin, 1995.

\bibitem{kamcss}
J.~Kamnitzer.
\newblock The crystal structure on the set of Mirkovi\'c-Vilonen polytopes.
\newblock {\tt arXiv:math.QA/0505398}.


\bibitem{katdcc}
J.~Kamnitzer and P.~Tingley.
\newblock A definition of the crystal commutor using {K}ashiwara's involution.
\newblock {\tt arXiv:math.QA/0610952}.

\bibitem{kascqa}
M.~Kashiwara.
\newblock Crystalizing the {$q$}-analogue of universal enveloping algebras.
\newblock {\em Commun. Math. Phys.}, 133:249--260, 1990.


\bibitem{kascb}
M.~Kashiwara.
\newblock On crystal bases.
\newblock In {\em Representations of groups (Banff, AB, 1994)}, volume~16 of
  {\em CMS Conf. Proc.}, pages 155--197. Amer. Math. Soc., Providence, RI,
  1995.

\bibitem{kancgr}
M.~Kashiwara and T.~Nakashima.
\newblock Crystal graphs for representations of the {$q$}-analogue of classical
  {L}ie algebras.
\newblock {\em J. Algebra}, 165:295--345, 1994.

\bibitem{lenccg}
C.~Lenart.
\newblock On the combinatorics of crystal graphs, {I}. {L}usztig's involution.
\newblock {\em Adv. Math.}, 211:204--243, 2007.

\bibitem{lapawg}
C.~Lenart and A.~Postnikov.
\newblock Affine {W}eyl groups in {$K$}-theory and representation theory.
\newblock {\tt arXiv:math.RT/0309207}.
\newblock To appear in {\em Int. Math. Res. Not.} 

\bibitem{lapcmc}
C.~Lenart and A.~Postnikov.
\newblock A combinatorial model for crystals of {K}ac-{M}oody algebras.
\newblock {\tt arXiv:math.RT/0502147}.
\newblock To appear in {\em Trans. Amer. Math. Soc.} 


\bibitem{litpro}
P.~Littelmann.
\newblock Paths and root operators in representation theory.
\newblock {\em Ann. of Math. {\rm (2)}}, 142:499--525, 1995.

\bibitem{litcrp}
P.~Littelmann.
\newblock Characters of representations and paths in {${\mathfrak h}\sp\ast\sb
  {\mathbb R}$}.
\newblock In {\em Representation theory and automorphic forms (Edinburgh,
  1996)}, volume~61 of {\em Proc. Sympos. Pure Math.}, pages 29--49. Amer.
  Math. Soc., Providence, RI, 1997.

\bibitem{luscba}
G.~Lusztig.
\newblock Canonical bases arising from quantized enveloping algebras. {II}.
\newblock {\em Progr. Theoret. Phys. Suppl.}, 102:175--201, 1991.

\bibitem{lusiqg}
G.~Lusztig.
\newblock {\em Introduction to Quantum Groups}, volume 110 of {\em Progress in
  Mathematics}.
\newblock Birkh\"auser Boston Inc., Boston, MA, 1993.

\bibitem{staec2}
R.~P. Stanley.
\newblock {\em Enumerative {C}ombinatorics. {V}ol. 2}, volume~62 of {\em
  Cambridge Studies in Advanced Mathematics}.
\newblock Cambridge University Press, Cambridge, 1999.

\bibitem{stecmw}
J.~R. Stembridge.
\newblock Combinatorial models for {W}eyl characters.
\newblock {\em Adv. Math.}, 168:96--131, 2002.

\bibitem{leuajt}
M.~A.~A. van Leeuwen.
\newblock An analogue of jeu de taquin for {L}ittelmann's crystal paths.
\newblock {\em S\'em. Lothar. Combin.}, 41:Art.\ B41b, 23 pp.\ (electronic),
  1998.

\end{thebibliography}

\end{document}